\documentclass[12pt]{article}
\usepackage[latin1]{inputenc}
\usepackage{amssymb,amsmath,amsfonts}
\newtheorem{theorem}{Theorem}
\newtheorem{proposition}[theorem]{Proposition}
\newtheorem{lemma}[theorem]{Lemma}

\newtheorem{corollary}[theorem]{Corollary}

\newcommand{\R}{\mathbb{R}}

\newcommand{\Sf}{\mathbb{S}}
\newcommand{\C}{\mathbb{C}}

\newcommand{\spa}{\mbox{span}}
\newcommand{\hess}{\mbox{Hess\,}}

\newcommand{\kerl}{\mbox{ker }}
\newcommand{\grad}{\mbox{grad\,}}

\newcommand{\po}{{\hspace*{-1ex}}{\bf .  }}

\newcommand{\nab}{\tilde\nabla}

\newcommand{\End}{\mbox{End}}
\newcommand{\trace}{\mbox{tr\,}}
\renewcommand{\o}{\omega}
\def\<{{\langle}}
\def\>{{\rangle}}

\def\T{{\cal T}}

\def\Y{{\cal Y}}
\def\n{\nabla}
\def\d{\partial}

\def\a{\alpha}

\def\bea{\begin{eqnarray*} }
\def\eea{\end{eqnarray*} }
\def\be{\begin{equation} }
\def\ee{\end{equation} }
\def\p{\partial}

\newcommand{\pu}{\partial_u}
\newcommand{\pv}{\partial_v}
\renewcommand{\gg}{\Gamma^1 }
\newcommand{\gh}{\Gamma^2 }
\def\proof{\noindent\emph{Proof: }}
\def\qed{\ifhmode\unskip\nobreak\fi\ifmmode\ifinner
\else\hskip5 pt \fi\fi\hbox{\hskip5 pt \vrule width4 pt
height6 pt  depth1.5 pt \hskip 1pt }}
\begin{document}

\title{The infinitesimally bendable Euclidean hypersurfaces}
\author{M.\ Dajczer and Th.\ Vlachos}
\date{}
\maketitle

\begin{abstract} 
The main purpose of this paper is to complete the work initiated 
by Sbrana in 1909  giving a complete local classification 
of the nonflat infinitesimally bendable hypersurfaces in Euclidean space.
\end{abstract}

 In the final decades of the 19th century geometers 
were increasingly interested in the study of  hypersurfaces in Euclidean 
space. Quite differently to what happens in the  surface 
case, these submanifolds are not easily isometrically deformable. 
In fact, it was shown that hypersurfaces are isometrically rigid 
provided that they bend in enough directions. 
The first correct proof that hypersurfaces with at least
three nonzero principal curvatures cannot be isometrically deformed was
given in 1885 by Killing \cite{ki} after a claim made in 1876 by 
Beez \cite{B}.

The situation of hypersurfaces of rank two, that is, the ones with exactly 
two nonzero principal curvatures, remained to be understood.  
It turned out that even in this case hypersurfaces are ``generically" rigid. 
After  earlier work for the three-dimensional case by Bianchi \cite{bi}, 
a parametric classification of all 
Euclidean  hypersurfaces $f\colon M^n\to\R^{n+1}$, $n\geq 3$, that admit 
non-trivial isometric deformations, was obtained in 1909 by Sbrana \cite{sb2}.  
This was done in terms of the so called Gauss parametrization discussed 
in the next section. Cartan \cite{ca} in 1916 gave a more careful statement  
but now in the language of envelopes of hyperplanes. 
See Dajczer, Florit and Tojeiro \cite{dft} for a modern presentation and
further results on the subject.

Perhaps the most interesting class in the classification discussed above 
is the one of isometrically bendable hypersurfaces, that is, when the 
hypersurface admits a smooth one-parameter variation by isometric hypersurfaces.  
These submanifolds can either be ruled, and then allow plenty of 
isometric bendings, or non-ruled in which case they just admit a single bending.

At around the same time,  Sbrana \cite{sb1} (who seems to have been a student
of Bianchi) in an  inspiring  paper considered the problem of classifying
Euclidean hypersurfaces that admit ``infinitesimal deformations", that is,
they are infinitesimally bendable. Roughly speaking, this  means
that there is a smooth one-parameter variation by hypersurfaces that are
isometric only to the ``first order". The precise definition of an 
infinitesimal bending is given in Section $3$. Of course, any bendable 
hypersurface is infinitesimally bendable, but the latter class turns out 
to be much larger. In fact, what Sbrana did was to provide a complete 
description of one class of infinitesimally bendable hypersurfaces 
(in terms of the Gauss parametrization already used in \cite{sb2}) but 
somehow ignored others.

It was very natural for Sbrana at that time to consider the infinitesimal 
version of the deformation problem. On one hand, because there was already 
a rich theory of infinitesimal bendings of surfaces; see Spivak \cite{sp}.  
On the other hand, it was known that any hypersurface
that possesses at least three nonzero principal curvatures is 
infinitesimally rigid, that is, it is not infinitesimally bendable,
a result that can be found in the book by Cesaro \cite{ce} from 1896.
A  modern proof of this fact follows from the more general result obtained 
by Dajczer and Rodr\'\i guez \cite{dr}.

It is for us quite surprising that we were not able to find any reference 
to Sbrana's contribution to the description of the hypersurfaces that admit 
infinitesimal bendings. In fact, the few places where his paper is referred to
are quite old and do not discuss his result; see \cite{sc} and \cite{st}.

We should point out that all of the above results are of local nature, as is
the case of this paper. By being local we mean that there is an open and dense
subset of the manifold such that along any connected component the submanifold
belongs to a class  in the classification. In that respect, we observe that
for isometric bendings  it was already shown  in \cite{dft} that hypersurfaces in
different classes can be smoothly attached.

The main purpose of this paper is to give a complete local classification of 
the nonflat infinitesimally bendable  hypersurfaces in modern terms.
In order to give a description of all hypersurfaces $f\colon M^n\to\R^{n+1}$,
$n\geq 3$, with two nonzero principal curvatures at any point that are 
infinitesimally bendable, we exclude from consideration the ones that are 
surface-like. Being \emph{surface-like} means that $f$ is locally part of
a cylinder either over a surface in $\R^3$ or the cone of a surface in 
$\Sf^3\subset\R^4$. The reason of exclusion is because in this case it can
be shown that the infinitesimal bending of the hypersurface is given by 
an infinitesimal bending of the surface, and the surface case is not an 
object of this paper.

Among the infinitesimally bendable  hypersurfaces there is the class of ruled 
hypersurfaces. A hypersurface $f\colon M^n\to\R^{n+1}$ is called \emph{ruled} 
if $M^n$ admits a foliation by leaves of codimension one mapped by $f$ into 
affine subspaces of $\R^{n+1}$. 
In our context, this class is not very interesting because it turns out 
that any infinitesimal bending is determined by an isometric bending.  
And isometric bendings of ruled hypersurfaces are easily seen 
to be parametrized by the set of smooth functions on an interval.

Finally, there is the class of infinitesimally bendable 
hypersurfaces that admit a unique infinitesimal bending. 
These hypersurfaces are the really interesting ones since generically they 
are not bendable, as we argue at the end of this introduction. 
We next give a characterization of the hypersurfaces belonging 
to this class \`{a} la Cartan, that is, in terms of envelopes of hyperplanes.  
An equivalent statement in terms of the Gauss parametrization, the one
used for the proof, is given later.
The concepts of envelope of hyperplanes and the Gauss
parametrization, as well as the relations between them, are be discussed 
in the next section. 

On an open subset $U\subset\R^2$ endowed with coordinates $(u,v)$ 
let  $\{\varphi_j\}_{0\leq j\leq n+1}$ 
be a set of solutions of the differential equation 
$$
\varphi_{z_1z_2}+M\varphi=0
$$
where $(z_1,z_2)$ can be either $(u,v)$ or $(u+iv,u-iv)$ and $M\in C^\infty(U)$. 
Assume that the map 
$\varphi=(\varphi_1,\ldots,\varphi_{n+1})\colon U\to\R^{n+1}$
is an immersion and consider the two-parameter family of affine hyperplanes 
$$
G(u,v)=\varphi_1x_1+\cdots +\varphi_{n+1}x_{n+1}-\varphi_0=0
$$
where $(x_1,\ldots,x_{n+1})$ are canonical coordinates of $\R^{n+1}$.

Our main result says that any hypersurface $f\colon M^n\to\R^{n+1}$ in 
the last class is the  envelope of a two-parameter family of hyperplanes 
as above which, in turn, means that $f$ is the solution of the system of 
equations
$
G=G_u=G_v=0.
$

\begin{theorem}\po\label{main0} Let $f\colon M^n\to\R^{n+1}$, $n\geq 3$, be an
infinitesimally bendable hypersurface of constant rank two that is neither
surface-like nor ruled on any open subset of $M^n$. Then, there is
an open and dense subset of $M^n$ such that along any connected
component $f$ is the envelope of a two-parameter family of hyperplanes as above.

Conversely, any hypersurface obtained as the envelope of a two-parameter family 
of hyperplanes as above admits locally a  unique infinitesimal bending.
\end{theorem}

Parametrically,  the hypersurface can be described by the  Gauss 
parametrization and goes as follows:
Let $g\colon U\to\Sf^n$ and $\gamma\in C^\infty(U)$ be given by
$$
g=\frac{1}{\|\varphi\|}(\varphi_1,\ldots,\varphi_{n+1})\;\;\mbox{and}
\;\;\gamma=\frac{\varphi_0}{\|\varphi\|}.
$$
If $\Lambda$ denotes the normal bundle of $g$ and  $h=i\circ g$ where 
$i\colon\Sf^n\to\R^{n+1}$  is the inclusion, then the map 
$\psi\colon \Lambda\to\R^{n+1}$ given by
$$
\psi(x,w)=\gamma(x)h(x)+h_*\grad\gamma(x)+w
$$
parametrizes the hypersurface.
\vspace{1ex}

We point out that for a hypersurface obtained as above, in order to be 
isometrically bendable the set of  functions $\varphi_1,\ldots,\varphi_{n+1}$ 
must satisfy a strong additional condition,
namely, the function $\phi=\|\varphi\|^2$
has to verify $\phi_{z_1z_2}=0$.

\section{Parametrizations}

In this section, we first recall how a Euclidean hypersurface of constant
rank can be locally parametrized by the use of the Gauss parametrization.
Then, we  discuss a class of envelopes of hyperplanes depending on parameters
as well how they can be  described in terms of the Gauss parametrization.

\subsection{The Gauss parametrization}

Let $f\colon M^n\to\R^{n+1}$ be an isometric immersion of constant rank 
$k$ for $1\leq k\leq n-1$.  By that we mean that its second fundamental 
form $A$ has constant rank $k$ or, equivalently, that the relative nullity
subspaces, i.e., the kernels of its second form  $\Delta(x)=\kerl A(x)$, 
satisfy $\dim\Delta(x)=n-k$  at any $x\in M^n$. 
In this situation, it is a standard fact that the tangent distribution 
$x\in M^n\mapsto\Delta(x)$ is integrable and that its totally geodesic 
leaves are mapped by $f$ into open subsets of affine subspaces of $\R^{n+1}$

A hypersurface of constant rank can be locally parametrized in terms of 
the image of its Gauss map $N$ and its support function $\gamma=\<f,N\>$.  
This parametrization is known as the  \emph{Gauss parametrization} 
and was described in \cite{dg} but it was already used by Sbrana  
in \cite{sb1} and \cite{sb2} long before.
\vspace{1ex}

Let  $(g,\gamma)$ be a pair formed by an isometric immersion 
$g\colon L^k\to\Sf^n$ into the unit sphere and a function 
$\gamma\in C^\infty(L)$. Denote by $\bar\pi\colon\Lambda\to L^k$ the normal 
bundle of $g$ and set $h=i\circ g$ where $i\colon\Sf^n\to\R^{n+1}$ 
is the standard inclusion. It was shown in \cite{dg} that the map 
$\psi\colon \Lambda\to\R^{n+1}$ given by
$$
\psi(x,w)=\gamma(x)h(x)+h_*\grad\gamma(x)+i_*w
$$
parametrizes (at regular points) a hypersurface of constant rank $k$ such
that the fibers of $\Lambda$ are identified with the leaves of the 
relative nullity foliation of $\psi$.

Conversely, any  hypersurface $f\colon M^n\to\R^{n+1}$ of constant  rank $k$ can be 
locally parametrized as above.   In fact, let  $U\subset M^n$ be 
an open saturated subset  of leaves of relative nullity
and let $\pi\colon U\to L^{n-k}$ denote the projection onto the quotient space.
The Gauss map $N$ of $f$  induces an immersion $g\colon L^{n-k}\to\Sf^n$ 
given by $g\circ\pi = N$. Moreover, since the support function  
$\<f,N\>$ is  constant along the relative nullity leaves, hence it induces 
a function $\gamma\in C^\infty(L)$. Now the Gauss parametrization allows 
to recover $f$ by means of the pair $(g,\gamma)$. 
\vspace{1,5ex}

The next statement presents some basic properties of the 
Gauss parametrization.

\begin{proposition}\po\label{gp} The following assertions hold:
\begin{itemize}
\item[(i)] The map $\psi\colon \Lambda\to\R^{n+1}$ is regular at $(x,w)$ 
if and only if  the self adjoint operator  
$$
P_w(x)=\gamma(x)I+\emph{Hess}\,\gamma(x)-A_w
$$ 
on $T_xL$ is nonsingular. Here $A_w$ is the shape operator of $g$ 
with respect to $w$. 
\item[(ii)]  The map $\psi$ when restricted to the open subset $V$ of regular 
points is an immersed hypersurface having the map $N\colon\Lambda\to\Sf^n$ 
given by $N(x,w)=g(x)$ as a Gauss map of rank $k$.
\item[(iii)] If $(x,w)\in V$ there is $j(x,w)\colon T_xL\to T_{(x,w)}\Lambda$ 
so that $j\colon T_xL\to\Delta^\perp(x,w)\subset T_{(x,v)}\Lambda$ is
an isometry such that
\be\label{eq}
h_*=\psi_*\circ j,\;\;\;P_w^{-1}=\bar\pi_*\circ j\;\;\mbox{and}\;\;
A\circ  j=-j\circ P_w^{-1}
\ee
where $A$ is the shape operator of $\psi$ at $(x,w)$ with respect to $N$.
\end{itemize}
\end{proposition}

\proof See \cite{dg}.\qed

\subsection{Envelopes of hyperplanes}

Let $P_u\colon\R^n\to\R^{n+1}$, $n\geq 2$, denote a smooth $k$-parameter 
family of affine hyperplanes parametrized by $u=(u_1,\ldots,u_k)$ on an 
open subset $U\subset\R^k$ with $1\leq k\leq n-1$. 
\vspace{1ex}

We say that a hypersurface  $f\colon M^n\to\R^{n+1}$ is the 
\emph{envelope of hyperplanes} of $(P_u)_{u\in U}$  if there exists 
a smooth totally geodesic foliation of $M^n$ by leaves 
$(L_u)_{u\in U}$ of dimension  $n-k$ parametrized by an embedding 
$h\colon U\to M^n$ transversal to the foliation and embeddings 
$j_u\colon L_u\to P_u(\R^n)$ as  open subset of $(n-k)$-dimensional 
affine subspaces of $\R^{n+1}$ such that $j_u=f\vert_{L_u}$ and 
$P_u(\R^n)=f_*T_{h(u)}M$ for any $u\in U$. 
\vspace{1ex}

Clearly, the leaves of $(L_u)_{u\in U}$ are contained in the relative nullity
of $f$. Also notice that any hypersurface $f\colon M^n\to\R^{n+1}$ with constant 
index of  relative nullity $\nu=n-k$ is the envelope of the $k$-parameter family of 
tangent hyperplanes.
\vspace{1ex} 

A $k$-parameter family of affine hyperplanes $(P_u)_{u\in U}$ can be given 
in terms of a smooth family  of equations of the form
$$
G(u)=\varphi_1x_1+\cdots +\varphi_{n+1}x_{n+1}-\varphi_0=0
$$
where $\varphi_j\in C^\infty(U)$, $0\leq j\leq n+1$, and  
$x=(x_1,\ldots,x_{n+1})$ are coordinates in $\R^{n+1}$ with respect 
to a canonical base.  

Assume that the map 
$\varphi=(\varphi_1,\ldots,\varphi_{n+1})\colon U\to\R^{n+1}$
is an immersion and without loss of generality that 
$0\not\in\varphi(U)$.  Let $g\colon U\to\Sf^n\subset\R^{n+1}$ 
and $\gamma\in C^\infty(U)$ be given by
$$
g=\frac{1}{\phi}(\varphi_1,\ldots,\varphi_{n+1})\;\;\mbox{and}
\;\;\gamma=\frac{\varphi_0}{\phi}
$$
where $\phi^2=\sum_{j=1}^{n+1}\varphi_j^2$. 

Being $g$ is an immersion, the pair $(g,\gamma)$  gives a 
hypersurface $f\colon M^n\to\R^{n+1}$ by means of the Gauss 
parametrization. Clearly $f$ is the envelope of $(P_u)_{u\in U}$ 
and the leaves $(L_u)_{u\in U}$ of the envelope coincide with the
relative nullity foliation of $f$.  Moreover,  the envelope of
$(P_u)_{u\in U}$ can be locally given as the solution of the 
system of equations
$$
(R)\begin{cases}
G(u)=0\\
G_{u_j}(u)=0,\;\;j=1,\ldots,k.
\end{cases}
$$  

We have shown the following fact.

\begin{proposition}\po Any hypersurface $f\colon M^n\to\R^{n+1}$
of constant rank $k$ can be locally  given as the envelope
of a smooth family of affine hyperplanes
$$
G(u)=\varphi_1x_1+\cdots +\varphi_{n+1}x_{n+1}-\varphi_0=0
$$
where $(\varphi_1,\ldots,\varphi_{n+1})\colon U\to\R^{n+1}$
is an immersion of an open subset $U$ of $\R^k$. 
Then $f$ is locally the solution of the  system of equations $(R)$.
\end{proposition}

\section{A class of surfaces}

A surface $g\colon L^2\to\Sf^n$ in the unit sphere is called
\emph{hyperbolic} (respectively, \emph{elliptic}) if there exists a
tensor $J$ on $L^2$ satisfying $J^2=I$ and $J\neq I$ 
(respectively, $J^2=-I$) and such that the second fundamental form 
$\a_g\colon TL\times TL\to N_gL$ of $g$ satisfies
\be\label{eq:aJ}
\a_g(JX,Y)=\a_g(X,JY)
\ee
for all vector fields $X,Y\in\mathfrak{X}(L)$. 
Local coordinates $(u,v)$ on $L^2$  are called
\emph{real-conjugate} for $g$ if the condition
$$
\alpha_g(\d_u,\d_v)=0
$$
holds where $\pu=\d/\d u$ and $\pv=\d/\d v$.
They are called \emph{complex-conjugate} if the condition
$\alpha_g(\p,\bar\p)=0$
holds where  $\p=\p_z=(1/2)(\pu-i\pv)$, that is, if we have that
$$
\alpha_g(\pu,\pu) + \alpha_g(\pv,\pv)=0.
$$

A simple argument (see \cite{dft}) gives the following result.

\begin{proposition}\po\label{coordinates}
Let $g\colon L^2\to\Sf^n$ be a hyperbolic (respectively, elliptic) surface.
Then  there exists locally a real-conjugate (respectively, complex-conjugate) 
system of coordinates on $L^2$ for $g$. Conversely, if there exists real-conjugate
(respectively, complex-conjugate) coordinates on $L^2$, then
$g$ is a hyperbolic (respectively, elliptic) surface.
\end{proposition}

Let $g\colon L^2\to\Sf^n$ be a simply-connected surface that carries a
real-conjugate system of coordinates $(u,v)$. Equivalently, the isometric 
immersion $h=i\circ g\colon L^2\to\R^{n+1}$ satisfies
\be\label{eq:guv}
h_{uv}-\gg h_u-\gh h_v + Fh=0
\ee
where $\gg,\gh$ are the Christoffel symbols given by
$$
\n_{\d_u}\d_v=\gg\d_u+\gh\d_v
$$
and $F=\<\d_u,\d_v\>$.

We are interested in surfaces for which, in addition, the following system of 
differential equations admits solution:
\be\label{eq:siti}
d\mu+2\mu\o=0\;\;\;\mbox{where}\;\;\; \o=\gh du+\gg dv.
\ee
This is the case if and only if the integrability condition 
\be\label{integcon}
\gg_u=\gh_v
\ee
is satisfied.

\begin{proposition}\po Let $g\colon L^2\to\Sf^n$ be a hyperbolic surface 
with real conjugate coordinates $(u,v)$ such that the induced metric satisfies
condition (\ref{integcon}).  Then there is a positive function 
$\mu\in C^\infty(L)$ such that $\varphi$ is a solution of (\ref{eq:guv})
if and only if $\psi=\sqrt{\mu}\varphi$ is a solution of 
$$
\psi_{uv}+M\psi=0
$$
where  $M\in C^\infty(L)$ is given by
\be\label{Mm}
M=F-\frac{\mu_{uv}}{2\mu}+\frac{\mu_u\mu_v}{4\mu^2}\cdot
\ee
In particular, the immersion 
$k=\sqrt{\mu}\,h\colon L^2\to\R^{n+1}$ satisfies
\be\label{form}
k_{uv}+Mk=0.
\ee
Conversely, let $k\colon L^2\to\R^{n+1}$ be an isometric immersion 
that for a system of coordinates  $(u,v)$ satisfies (\ref{form})
where  $M\in C^\infty(L)$. Then $(u,v)$ are real conjugate coordinates
for the immersion $g=(1/\|k\|)\,k\colon L^2\to\Sf^n$ and
condition (\ref{integcon}) is satisfied for the induced metric.
\end{proposition}

\proof If $g$ satisfies the integrability condition, then
$$
\gg=-\frac{\mu_v}{2\mu},\;\;\;\gh=-\frac{\mu_u}{2\mu}
$$
where $\mu=c\,e^{-2\int\o}$ for any $c\in\R_+$. Hence (\ref{eq:guv}) 
becomes
\be\label{eq:guv2}
h_{uv}+\frac{\mu_v}{2\mu}h_u+\frac{\mu_u}{2\mu}h_v+Fh=0.
\ee
It follows easily that $k=\sqrt{\mu}\,h$ takes the form (\ref{form}) 
where $M$ is given by (\ref{Mm}).
The converse is a straightforward computation.\vspace{1,5ex}\qed

Let $g\colon L^2\to\Sf^n$ be a simply-connected surface endowed
with  complex-conjugate coordinates $(z,\bar{z})$.  
Equivalently, the isometric immersion $h=i\circ g\colon L^2\to\R^{n+1}$
satisfies
\be\label{eq:guz}
h_{z\bar z}-\Gamma h_z-\bar\Gamma h_{\bar z}+Fh=0
\ee
where the Christoffel symbols, obtained using the $\C$-linear extensions 
of the metric of $L^2$ and the corresponding connection, are given by
$$
\n_{\d}\bar\d=\Gamma\d+\bar \Gamma\bar\d
$$ 
and $F=\<\d,\bar\d\>$.

We are interested in surfaces for which, in addition, the following 
system of differential equations for $\mu$ real admits solutions:
\be \label{eq:comp}
\mu_{\bar{z}}+2\mu\Gamma=0.
\ee
This is the case if and only if the integrability condition 
\be\label{integcon2}
\Gamma_z=\bar{\Gamma}_{\bar z}, 
\ee
that is, $\Gamma_z\;\mbox{is real}$, is satisfied.

\begin{proposition}\po Let $g\colon L^2\to\Sf^n$ be an elliptic surface
with complex  conjugate coordinates $(z,\bar{z})$  such that the induced 
metric satisfies condition (\ref{integcon2}). 
Then, there is a positive solution  $\mu\in C^\infty(L)$ 
of (\ref{eq:comp}) such that $\varphi$ is a solution of (\ref{eq:guz})
if and only if $\psi=\sqrt{\mu}\varphi$ is a solution of 
$$
\psi_{z\bar z}+M\psi=0
$$
where  $M\in C^\infty(L)$ is given by
\be\label{Mm2}
M=F-\frac{\mu_{z\bar z}}{2\mu}+\frac{\mu_z\mu_{\bar z}}{4\mu^2}\cdot
\ee
In particular, the immersion 
$k=\sqrt{\mu}\,h\colon L^2\to\R^{n+1}$ satisfies
\be\label{form2}
k_{z\bar z}+Mk=0.
\ee

Conversely, let $k\colon L^2\to\R^{n+1}$ be an isometric immersion that 
 for a system of coordinates  $(z,\bar{z})$ satisfies (\ref{form2})
where  $M\in C^\infty(L)$. Then $(z,\bar{z})$ are complex conjugate
coordinates for the immersion $g=(1/\|k\|)k\colon L^2\to\Sf^n$  and 
condition  (\ref{integcon2}) is satisfied for the induced metric.
\end{proposition}

\proof We have $\mu=c\,e^{-2\int\o}$ for any $c\in\R_+$ where
$\o=\Gamma d\bar{z}$. Then (\ref{eq:guz}) takes the form
$$
h_{z\bar z}+\frac{\mu_z}{2\mu}h_z+\frac{\mu_{\bar z}}{2\mu}h_{\bar z}+Fh=0.
$$
It follows easily that $k=\sqrt{\mu}\,h$ is as in  (\ref{form2}) where 
$M$ is given by (\ref{Mm2}).
The converse is a straightforward computation.\qed

\section{The main result}

After introducing the necessary terminology and definitions, we  present 
the main result of the paper in terms of the Gauss parametrization, as 
is the case in the paper by Sbrana.  The proof of the alternative version 
of the theorem in terms  of envelopes of hyperplanes given in the introduction 
can easily be obtained from this version using results from the preceding 
sections. 
\vspace{1,5ex}

By a \emph{variation} $F$ of an isometric immersion $f\colon M^n\to\R^{n+1}$ 
we mean a smooth map $F\colon (-\epsilon, \epsilon)\times M^n\to\R^{n+1}$ 
such that $f_t=F(t,\cdot)$ is an immersion for each $t\in I=(-\epsilon,
\epsilon)$ and  $f=f_0$.
The variational vector field of $F$ is the section $\T\in\Gamma(f^*(T\R^{n+1}))$
of the Riemannian vector bundle $f^*(T\R^{n+1})$ defined as
$$
\T(x)=F_*\d/\d t|_{t=0}(x).
$$

A variation $F$ of a given isometric immersion
$f\colon M^n\to\R^{n+1}$ is called an \emph{isometric bending} if 
$f_t$ is an isometric immersion for any $t\in I$.
The variational vector field of an isometric bending satisfies
$$
\<\nab_X\T,f_*Y\>+\<f_*X,\nab_Y\T\>=0
$$
for any $X,Y\in\mathfrak{X}(M)$. Equivalently, it satisfies that
$$
\<\nab_X\T,f_*X\>=0
$$
for any $X\in\mathfrak{X}(M)$.
\vspace{1ex}

An isometric bending $F$ is called \emph{trivial} if it is produced by a smooth
one-parameter family of isometries of $\R^{n+1}$, that is, if there
exist a smooth family ${\cal C}\colon I\to O(n+1)$ of orthogonal transformations
of $\R^{n+1}$ and a smooth map $v\colon I\to \R^{n+1}$ such that
$$
F(t,x)={\cal C}(t)f(x)+v(t).
$$
For a trivial isometric bending the variational vector field is of the form
$$
\T(x)={\cal D}f(x)+w
$$
where ${\cal D}={\cal C}'(0)$ is a skew-symmetric linear endomorphism of $\R^{n+1}$
and $w=v'(0)$ a vector in  $\R^{n+1}$.
Conversely, given a skew-symmetric linear endomorphism ${\cal D}$ of $\R^{n+1}$
and a vector $w\in\R^{n+1}$, the map
$$
F(t,x)=e^{t{\cal D}}f(x)+tw
$$
defines a  trivial isometric bending that has $\T={\cal D}f+w$ as variational
vector field.
\vspace{1ex}

By an \emph{infinitesimal bending} $\T$ of an isometric immersion
$f\colon M^n\to\R^{n+1}$ we mean an element of $\Gamma(f^*(T\R^{n+1}))$
that satisfies
\be\label{iif0}
\<\nab_X\T,f_*Y\>+\<f_*X,\nab_Y\T\>=0
\ee
for any $X,Y\in\mathfrak{X}(M)$.
An infinitesimal bending is said to be \emph{trivial} if
$$
\T(x)={\cal D}f(x)+w
$$
where ${\cal D}$ is a skew-symmetric linear
endomorphism of $\R^{n+1}$ and $w\in\R^{n+1}$.
\vspace{1ex}

An isometric immersion $f\colon M^n\to\R^{n+1}$ is called \emph{infinitesimally
bendable} if it admits a nontrivial infinitesimal bending.  Otherwise, 
it is said that $f$ is \emph{infinitesimally rigid}.
\vspace{1ex}

Multiplying a given infinitesimal bending by a real constant and adding a trivial
infinitesimal bending yields a new infinitesimal bending. In the sequel,
we identify two infinitesimal bendings $\T_1$ and $\T_2$ if
$\T_2=\T_0 + c\,\T_1$ where $\T_0$ is a trivial infinitesimal bending and 
$0\neq c\in\R$.
\vspace{1ex}

We have already observed  that hypersurfaces of rank at least three at any point 
are infinitesimally rigid.  Therefore, the interesting case to be considered is 
the one of constant rank two.  We see next that even in this special case  
hypersurfaces are ``generically" infinitesimally rigid.
\vspace{1ex}

We call the pair $(g,\gamma)$ a \emph{special hyperbolic pair}
(respectively, \emph{special elliptic pair}) if  $g\colon L^2\to\Sf^n$
is a hyperbolic (respectively, elliptic) surface so that system (\ref{eq:siti})
(respectively, system (\ref{eq:comp})) has solution and $\gamma\in C^\infty(L)$ 
satisfies (\ref{eq:guv}) (respectively,  (\ref{eq:guz})).

\begin{theorem}\po\label{main} Let $f\colon M^n\to\R^{n+1}$, $n\geq 3$, be an
infinitesimally bendable hypersurface of constant rank two that is neither
surface-like nor ruled on any open subset of $M^n$. Then, there is
an open and dense subset of $M^n$ such that along any connected
component $f$ is parametrized in terms of the Gauss parametrization
by a special hyperbolic or a special elliptic pair.

Conversely, any hypersurface parametrized in terms of the Gauss
parametrization by a special hyperbolic or special elliptic pair
admits locally a  unique infinitesimal bending.
\end{theorem}

The case of ruled hypersurfaces that has been excluded from consideration 
in the above result is rather simple and will be treated separately in 
Section $5$.

\section{Existence and uniqueness}

 We study the system of differential equations of  an infinitesimal 
bending of a Euclidean hypersurface and discuss its integrability
conditions. This yields a kind of fundamental theorem for 
infinitesimal bendings that is,  basically, contained in Sbrana's 
paper  \cite{sb1}.  In fact, the case of arbitrary codimension 
was later taken on by Schouten \cite{sc} but presented in a rather difficult 
terminology.  We point out that in this section some  long but straightforward 
computations are only indicated.
\vspace{1,5ex}

Given a hypersurface $f\colon M^n\to\R^{n+1}$, $n\geq 3$, in the sequel 
we associate to any infinitesimal bending $\T$ of $f$ the variation
$F\colon \R\times M^n\to\R^{n+1}$ with variational vector field 
$\T$ given by
$$
F(t,x)=f(x)+t\T(x).
$$
It is usually said that  $f_t=F(t,\cdot)$ is isometric to $f$ up to first
order for if
\be\label{var}
\|f_{t*}X\|^2=\|f_*X\|^2+t^2\|\nab_X\T\|^2
\ee
for all $X\in\mathfrak{X}(M)$.

Let $g_t$ be the metric on $M^n$ induce by $f_t$. Then,
$$
\p/\p t\vert_{t=0}\,g_t(X,Y)=0
$$
for all $X,Y\in\mathfrak{X}(M)$.  Consequently, we have that the  associated
one-parameter family of Levi-Civita connections and the corresponding
family of curvature tensors satisfy 
$$
\p/\p t\vert_{t=0}\,f_*\n^t_XY=0
$$
and
\be\label{curvature}
\p/\p t\vert_{t=0}\,\<R^t(X,Y)Z,W\>=0
\ee
for all $X,Y,Z,W\in\mathfrak{X}(M)$.
\vspace{1ex}

 Let $N(t)$ denote a Gauss map of $f_t$ and $A(t)$ the second fundamental 
form of $f_t$ with respect to $N(t)$ so that the map $t\in\R\mapsto N(t)$ 
is smooth.  Then  $N=N(0)$ is the Gauss map and $A=A(0)$ is the second 
fundamental form of $f$. Moreover, let us define
$L\in\Gamma(\End(TM,f^*(T\R^{n+1}))$ by
$$
LX=\nab_X\T=\T_*X.
$$
Then (\ref{iif0}) can be written as
\be\label{iif}
\<LX,f_*Y\>+\<f_*X,LY\>=0
\ee
for all $X,Y\in\mathfrak{X}(M)$.

\begin{lemma}\po We have that 
$\Y=\p/\p t\vert_{t=0}\, N(t)\in\Gamma(f^*(T\R^{n+1}))$ satisfies
\be\label{4}
\<\Y,N\>=0
\ee
and
\be\label{3}
\<\Y,f_*X\>+\<LX,N\>=0
\ee
for all $X\in\mathfrak{X}(M)$.
\end{lemma}

\proof The derivative with respect to $t$ at $t=0$ of $\<N(t),N(t)\>=1$ 
gives (\ref{4}) whereas  of $\<N(t),f_{t*}X\>=0$ yields (\ref{3}).\qed

\begin{lemma}\po  
We have that $B=\p/\p t\vert_{t=0}\, A(t)\in\Gamma(\End(TM))$ 
is symmetric and satisfies 
\be\label{1}
(\nab_XL)Y=\<BX,Y\>N+\<AX,Y\>\Y
\ee
and
\be\label{2}
\Y_*X=-f_*BX-LAX
\ee
for all $X,Y\in\mathfrak{X}(M)$.
\end{lemma}

\proof The derivative with respect to $t$ at $t=0$ of the Gauss formula
$$
\nab_Xf_{t*}Y=f_{t*}\n^t_XY + g_t(A(t)X,Y)N(t)
$$
 easily gives (\ref{1}).
As for the Weingarten formula
$$
\nab_XN(t)=-f_{t*}(A(t)X)
$$
we have that its derivative at $t=0$ yields (\ref{2}).
\vspace{1,5ex}\qed

If $\T={\cal D}f+w$ is a trivial infinitesimal bending then 
$L={\cal D}\circ f_*$.  It follows that 
$\Y={\cal D}N$  and that $B=0$ since 
$$
\<BX,Y\>=\<(\nab_XL)Y,N\>=\<(\nab_X{\cal D})Y,N\>=0
$$
for all $X,Y\in\mathfrak{X}(M)$.

\begin{proposition}\po  The tensor $B$ is a symmetric Codazzi 
tensor, i.e.,
\be\label{6}
\left(\n_XB\right)Y-\left(\n_YB\right)X=0
\ee
such that
\be\label{5}
BX\wedge AY-BY\wedge AX=0
\ee
for all $X,Y\in\mathfrak{X}(M)$.
\end{proposition}

\proof  The derivative at $t=0$ of the Codazzi equation
$$
\left(\n^t_XA(t)\right)Y=\left(\n^t_YA(t)\right)X
$$
gives (\ref{6}).
To obtain (\ref{5}) we compute the derivative at $t=0$ of the
Gauss equation
$$
R^t(X,Y)Z=g(t)(A(t)Y,Z)A(t)X-g(t)(A(t)X,Z)A(t)Y
$$
and use (\ref{curvature}).\vspace{1,5ex}\qed

The next result is to be expected bearing in mind the nature of the Gauss
and Codazzi equations as the integrability conditions for the system of
differential equations associated to an isometric immersion as a hypersurface.

\begin{lemma}\po\label{sys} Equations (\ref{6}) and (\ref{5}) are the
integrability conditions of the system of  differential equations
(\ref{1}) and (\ref{2}) for $L$ and $\Y$, that is,
$$
(S)\begin{cases}
\Y_*X=-LAX-f_*BX\\
(\nab_XL)Y=\<BX,Y\>N+\<AX,Y\>\Y.
\end{cases}
$$
\end{lemma}

\proof For the first equation, we have to show that
\be\label{condi1}
\nab_X\Y_*Y-\nab_Y\Y_*X-\nab_{[X,Y]}\Y=0
\ee
for all $X,Y\in\mathfrak{X}(M)$. One has that
$$
\nab_X\Y_*Y=-(\nab_XL)AY-L(\n_XA)Y-LA\n_XY-f_*\n_XBY-\<AX,BY\>N.
$$
Then (\ref{condi1}) is equivalent to
$$
(\nab_XL)AY-(\nab_YL)AX+f_*((\n_XB)Y-(\n_YB)X)
+(\<AX,BY\>-\<AY,BX\>)N=0.
$$
Replacing the first two terms by the use of the second equation in $(S)$ 
it is easily seen that (\ref{condi1}) follows from (\ref{6}).

It is easy to see that the integrability condition for the second
equation is
\be\label{condi2}
(\nab_X\nab_YL-\nab_Y\nab_XL-\nab_{[X,Y]}L)Z=-LR(X,Y)Z
\ee
for all $X,Y,Z,W\in\mathfrak{X}(M)$.
A straightforward computation using (\ref{1}) gives
\bea
(\nab_X\nab_YL)Z\!\!\!&=&\!\!\!\<(\n_XB)Y,Z\>N+\<B\n_XY,Z\>N
-\<BY,Z\>f_*AX+\<(\n_XA)Y,Z\>\Y\\
\!\!\!&&\!\!\!+\; \<A\n_XY,Z\>\Y
-\<AY,Z\>LAX-\<AY,Z\>f_*BX.
\eea
That $A$ is a Codazzi tensor together with (\ref{6}) yields
\bea
(\nab_X\nab_YL-\nab_Y\nab_XL-\nab_{[X,Y]}L)Z
\!\!\!&=&\!\!\!-\<BY,Z\>f_*AX-\<AY,Z\>(LAX+f_*BX)\\
\!\!\!&&\!\!\!+\,\<BX,Z\>f_*AY-\<AX,Z\>(LAY+f_*BX).
\eea
On the other hand,  we have
$$
LR(X,Y)Z=\<AY,Z\>LAX-\<AX,Z\>LAY,
$$
and (\ref{condi2}) follows using (\ref{5}).\vspace{1,5ex}\qed

Next we consider the case of hypersurfaces of constant rank two.

\begin{corollary}\label{kerB}\po If $f\colon M^n\to\R^{n+1}$, $n\geq 3$, 
is an infinitesimally bendable hypersurface of constant rank two, 
then $\Delta\subset\ker B$.
\end{corollary}

\proof  This follows easily from (\ref{5}).

\begin{theorem}\po\label{mains}  Let $f\colon M^n\to\R^{n+1}$, $n\geq 3$, be a 
simply-connected  hypersurface of constant rank two.  
Then, the set of all symmetric Codazzi tensors
$B\in\Gamma(\End(TM))$  such that $\Delta\subset\ker B$ and
$$
BX\wedge AY-BY\wedge AX=0
$$
for all $X,Y\in\mathfrak{X}(M)$ is in one-to-one correspondence with 
the set of all infinitesimal bendings of $f$ so that $B=0$ corresponds 
to the trivial one.
\end{theorem}

\proof Given  $B\in\Gamma(\End(TM))$ as in the statement,  
we first prove that there exists a solution $\Y$  and $L$ of system 
$(S)$ such that (\ref{iif}), (\ref{4}) and (\ref{3}) are satisfied. 
In particular, this gives the existence of an infinitesimal bending 
$\T$ such that $L=\T_*$. To see this, observe that by (\ref{1}) the 
one-form $\o=\<L,v\>$ is closed for any $v\in\R^{n+1}$.

Given a solution  $\Y$ and $L$ of $(S)$, we define a smooth function by
$$
\tau=\<\Y,N\>,
$$
a smooth one-form by
$$
\theta(X)=\<\Y,f_*X\>+\<LX,N\>
$$
and a smooth symmetric bilinear tensor by
$$
\beta(X,Y)=\<LX,f_*Y\>+\<LY,f_*X\>.
$$
A straightforward calculation gives that
\be\label{tau}
d\tau=-\theta\circ A,
\ee
\be\label{theta}
(\n_X\theta)Y=-\beta(AX,Y)+2\tau\<AX,Y\>
\ee
and
\be\label{M}
(\n_Z\beta)(X,Y)=\<AX,Y\>\theta(Z)+\<AX,Z\>\theta(Y)
\ee
for all $X,Y,Z\in\mathfrak{X}(M)$.

We claim that the system of differential equations formed by (\ref{tau}), 
(\ref{theta}) and (\ref{M}) is completely integrable.  The integrability 
condition for the first equation is easy to verify. For the second equation, 
we have to see that
\be\label{cal1}
(\n_X\n_Y\theta-\n_Y\n_X\theta-\n_{[X,Y]}\theta)Z
=-\theta(R(X,Y)Z)
\ee
holds. Using (\ref{tau}) and (\ref{theta}) we obtain 
$$
(\n_X\n_Y\theta)Z=-(\n_X\beta)(AY,Z)
-\beta(\n_XAY,Z)-2\theta(AX)\<AY,Z\>
+2\tau\<\n_XAY,Z\>.
$$
Hence
\begin{align*}
 (\n_X\n_Y\theta-\n_Y\n_X\theta-\n_{[X,Y]}\theta)Z
=&-(\n_X\beta)(AY,Z)+(\n_Y\beta)(AX,Z)\\
&\,-2\theta(AX)\<AY,Z\>+2\theta(AY)\<AX,Z\>.
\end{align*}
Using (\ref{M}) we obtain that
$$
(\n_X\n_Y\theta-\n_Y\n_X\theta-\n_{[X,Y]}\theta)Z
=-\theta(AX)\<AY,Z\>+\theta(AY)\<AX,Z\>.
$$
On the other hand, we have from the Gauss equation that
$$
\theta(R(X,Y)Z)=\<AY,Z\>\theta(AX)-\<AX,Z\>\theta(AY),
$$
and (\ref{cal1}) follows.

Finally, the integrability condition for the last equation, namely, that  
\be\label{cal2}
\!\!(\n_X\n_Y\beta-\n_Y\n_X\beta-\n_{[X,Y]}\beta)(Z,W)
=-\beta(R(X,Y)Z,W)-\beta(R(X,Y)W,Z)
\ee
can be verified by a similar computation, and this
proves the claim.

Start with a solution $L^*$ and $\Y^*$ of system $(S)$ with 
corresponding tensors $\theta^*$,  $\beta^*$ and function $\tau^*$.
Fix a point $p_0\in M^n$ and let $L_0$ and $\Y_0$ be a solution of the 
integrable system
$$
(S_0)\begin{cases}
\Y_*X=-LAX\\
(\nab_XL)Y=\<AX,Y\>\Y
\end{cases}
$$
with  initial conditions $\theta_0(p_0)=\theta^*(p_0)$, 
$\beta_0(p_0)=\beta^*(p_0)$ and $\tau_0(p_0)=\tau^*(p_0)$.
Then $L=L^*-L_0$ and
$\Y=\Y^*-\Y_0$ are a solution of $(S)$ such that 
$\theta=\theta^*-\theta_0$,  $\beta=\beta^*-\beta_0$ and 
$\tau=\tau^*-\tau_0$. Clearly $\theta(p_0)=\beta(p_0)=\tau(p_0)=0$.
Since $\theta, \beta$ and $\tau$ solve the homogeneous integrable system
(\ref{tau}), (\ref{theta}) and (\ref{M}), hence $\theta=\beta=\tau=0$.

Given any two pairs  $L_j,\Y_j$, obtained as above, let $\T_j$, $1\leq j\leq 2$, 
be the associated infinitesimal bendings. 
It remains to show that $\T=\T_1-\T_2$ is a trivial infinitesimal bending.

We have that the pair $L=L_1-L_2$, $\Y=\Y_1-\Y_2$ satisfies $(S_0)$
as well as  (\ref{iif}), (\ref{4}) and (\ref{3}). Fix $p_0\in M^n$ and
define  a skew-symmetric linear endomorphism ${\cal C}$ of $\R^{n+1}$ by
$$
{\cal C}f_*(p_0)X=L(p_0)X\;\;\mbox{and}\;\;{\cal C}N(p_0)=\Y(p_0)
$$
and a vector $v\in\R^{n+1}$ by $v=\T(p_0)-{\cal C}f(p_0)$.
Consider the  trivial infinitesimal bending $\tilde{\T}={\cal C}f+v$ and
$\tilde{\Y}={\cal C}N$.  Then, the pair $\tilde{L}$ and $\tilde{\Y}$ satisfies
$(S_0)$. Thus, also the pair $L^*=L-\tilde{L}$, $\Y^*=\Y-\tilde{\Y}$
solves system $(S_0)$. Moreover, $\T^*(p_0)=0$, $\Y^*(p_0)=0$ and
$L^*(p_0)=L(p_0)-\tilde{L}(p_0)=0$.
Thus $\T^*=0$ and hence $\T=\tilde{\T}$.\qed

\section{The proof of Theorem \ref{main}}

In the sequel, let $f\colon M^n\to\R^{n+1}$, $n\geq 3$, be a 
hypersurface of constant  rank two. Recall that the \emph{splitting tensor} 
$C\colon\Gamma(\Delta)\to\Gamma(\End(\Delta^\perp))$
is defined by 
$$
C_TX=-(\n_XT)_{\Delta^\perp}
$$
for any $T\in\Gamma(\Delta)$ and $X\in\mathfrak{X}(M)$. 
 From the Codazzi equation, it follows that
\be\label{codazzi}
\n_TA=AC_T=C_T^tA
\ee
for any $T\in\Gamma(\Delta)$.

\begin{proposition}\po  Assume that the splitting tensor at any point
satisfies $C_T\in\spa\{I\}$ for any $T\in\Delta$, where $I$ denotes the 
identity section of $\End(\Delta^\perp)$. Then $f$ is surface-like.
\end{proposition}

\proof See Lemma $6$ in \cite{dft}.\vspace{1,5ex}\qed

Assume further that $f$ is infinitesimally bendable. 
Locally and because of the rank assumption,  there is an orthonormal
tangent frame spanning $\Delta^\perp$ such that
\be\label{matrix}
A\vert_{\Delta^\perp}=\begin{bmatrix}
\lambda_1&0\\
0&\lambda_2
\end{bmatrix}.
\ee

\begin{lemma}\po If $B\neq 0$ at any point of $M^n$, then 
\be\label{matrix2}
B\vert_{\Delta^\perp}=\begin{bmatrix}
\lambda\lambda_1&b\\
b&-\lambda\lambda_2
\end{bmatrix}.
\ee
\end{lemma}

\proof By Corollary \ref{kerB} we have that $\Delta\subset\ker B$.
Now (\ref{matrix2}) follows easily from (\ref{5}).\qed

\begin{lemma}\po\label{conditions} We have that 
$D=(A\vert_{\Delta^\perp})^{-1}B\vert_{\Delta^\perp}\in\Gamma(\End(\Delta^\perp))$ 
satisfies:
\begin{itemize}
\item[(i)] $[D,C_T]=0$ for all $T\in\Delta$,
\item[(ii)] $\n_TD=0$ for all $T\in\Delta$,
\item[(iii)] $\trace D=0$,
\item[(iv)] $T(\det D)=0$ for all $T\in\Delta$.
\end{itemize}
\end{lemma}

\proof We denote $A=A\vert_{\Delta^\perp}$ and $B=B\vert_{\Delta^\perp}$.   
From (\ref{codazzi}) we obtain $\n_TB=BC_T.$
Hence
$$
BC_T=C_T^tB.
$$
We have using  (\ref{codazzi}) that
$$
ADC_T=BC_T=C_T^tB=C_T^tAD=AC_TD,
$$
and $(i)$ follows. We have
\bea
A\n_TD\!\!\!&=&\!\!\!\n_T(AD)-(\n_TA)D=\n_TB-(\n_TA)D
=BC_T-AC_TD\\
\!\!\!&=&\!\!\!BC_T-C_T^tAD=BC_T-C_T^tB=0,
\eea
and this yields $(ii)$. We obtain from (\ref{matrix}) and (\ref{matrix2}) that
$$
D=\begin{bmatrix}
\lambda&b/\lambda_1\\
b/\lambda_2&-\lambda
\end{bmatrix},
$$
which gives $(iii)$. Now part $(iv)$ follows
from  $(ii)$ and $(iii)$.\qed

\begin{proposition}\po\label{ruled} Assume that $f\colon M^n\to\R^{n+1}$, $n\geq 3$,
is not surface-like on any open subset of $M^n$. Then $f$ is ruled
along any open subset where $D\neq 0$ satisfies $\det D=0$.
\end{proposition}

\proof By Lemma \ref{conditions} there is an orthogonal frame $X,Y$
of $\Delta^\perp$ with $Y$ of unit length such that $DY=0$ and $DX=Y$.
We claim that $f$ is ruled by the integral leaves of
the distribution
$\Delta\oplus\spa\{Y\}$. To see this, we have to show that
$$
(i)\; \<AY,Y\>=0,\;(ii)\;\n_TY=0,\;
(iii)\;\<\n_YT,X\>=0\;\;\mbox{and}\;\;(iv)\;\<\n_YY,X\>=0
$$
for all $T\in\Delta$. We have
$$
\<AY,Y\>=\<ADX,Y\>=\<BX,Y\>=\<BY,X\>=\<ADY,X\>=0.
$$
Condition $(ii)$ follows easily using $\n_TD=0$.
Since  $[D,C_T]=0$, we obtain
$$
\<\n_YT,X\>=-\<C_TY,X\>=-\<C_TDX,X\>=-\<DC_TX,X\>=0.
$$
We have that
$$
BY=ADY=0\;\;\mbox{and}\;\;BX=ADX=AY=\lambda X,\;\;\lambda\neq 0,
$$
and condition $(iv)$ follows easily using (\ref{6}).
\vspace{1,5ex}\qed

In the sequel, we consider the case $\det D\neq 0$. By the above, this 
is always the case under the assumptions of Theorem \ref{main}.  
\vspace{1ex}

By part $(iii)$ of Lemma \ref{conditions} the eigenvalues of $D$ are 
the solutions of  $t^2+\det D=0$. Therefore, on each connected component 
of an open subset of $ M^n$ either $\det D<0$ and thus $D$ has two 
smooth real eigenvalues $\{\mu,-\mu\}$ or $\det D>0$ and thus $D$ has a 
pair of smooth complex eigenvalues $\{i\mu,-i\mu\}$. Then 
$J\in\Gamma(\End(\Delta^\perp))$ defined by 
\be\label{jota}
D=\mu J
\ee
satisfies $J^2=I$ in the first case and $J^2=-I$ in the second case.

\begin{lemma}\po\label{para}  The eigenspaces of $D$ are parallel and the
eigenvalues constant along the leaves of $\Delta$.
\end{lemma}

\proof Follows from parts $(ii)$ and $(iv)$ of Lemma \ref{conditions}.
\vspace{1,5ex}\qed

A hypersurface $f\colon M^n\to\R^{n+1}$ of rank two is said to be
\emph{hyperbolic} (respectively, \emph{elliptic}) if there exists
$J\in\Gamma(\End(\Delta^\perp))$ satisfying the following conditions:
\begin{itemize}
\item[(i)]  $J^2=I$ and $J\neq I$ (respectively, $J^2=-I$).
\item[(ii)] $\n_TJ=0$ for all $T\in\Gamma(\Delta)$.
\item[(iii)] $C_T\in\spa\{I,J\}$ for all $T\in\Gamma(\Delta)$.
\end{itemize}

\begin{proposition}\label{hypel}
\po Assume that $f\colon M^n\to\R^{n+1}$ is neither surface-like
nor ruled on any open subset $\tilde{M}^n$ of $M^n$. Then, there is an 
open and dense subset $\tilde M^n$ of $M^n$ such that the restriction of 
$f$ to any connected component of $\tilde{M}^n$ is either hyperbolic or 
elliptic.
\end{proposition}

\proof Let  $J\in\Gamma(\End(\Delta^\perp))$ be defined by (\ref{jota}).
The subspace $S$ of all elements in $\End(\Delta^\perp)$ that
commute  with $D$, i.e., that commute with $J$, is $S=\spa\{I,J\}$.
Thus condition $(iii)$ in the above definition follows from part 
$(i)$ of Lemma \ref{conditions}.\vspace{1,5ex}\qed

Given a submersion $\pi\colon M\to L$ between differentiable manifolds,
then  $X\in \mathfrak{X}(M)$ is said  to be \emph{projectable} if it is 
$\pi$-related to some  $\bar X\in \mathfrak{X}(L)$, that is, if there 
exists $\bar X\in \mathfrak{X}(L)$ such that $\pi_*X=\bar X\circ \pi$.
\vspace{1,5ex}

In the sequel, we denote by $\pi\colon M^n\to L^2$ the submersion onto the
(local) quotient space of leaves of $\Delta$, namely, onto $L^2=M^n/\Delta$.
A tensor $D\in\End(\Delta^\perp)$  is said to be  \emph{projectable} with
respect to $\pi$ if it is the horizontal lift of some tensor $\bar D$ on $L$.
Clearly, $D$ is projectable with respect to  $\pi$ if and only if for all
$\bar{x}\in L$,  $x,y\in \pi^{-1}(\bar{x})$,
$v\in\Delta^\perp(x)$ and $w\in \Delta^\perp(y)$ with  $\pi_*v=\pi_*w$,
we have  that $\pi_*Dv=\pi_*Dw$.

\begin{lemma}\label{prop*}\po
Let $f\colon M^n\to\R^{n+1}$ be a  hypersurface of rank two parametrized  
by a pair $(g,\gamma)$ in terms of the  Gauss parametrization.
If $f$ is  hyperbolic (respectively, elliptic) with respect to  
$J\in\Gamma(\End(\Delta^\perp))$ and $D=\mu J$ satisfies  $(i)$--$(iv)$ 
in Lemma~\ref{conditions}, then $J$ and $D$ are the horizontal lifts of 
tensors $\bar{J}$ and $\bar D=\bar\mu\bar{J}$ on $L^2$ such that 
$\mu=\bar{\mu}\circ\pi$, $\bar{J}^2=\bar I$ 
(respectively, $\bar{J}^2=-\bar I$), the pair  $(g,\gamma)$ is hyperbolic  
(respectively, elliptic)  with respect to  $\bar{J}$ and $\bar D$ satisfies:
\begin{itemize}
\item[(a)] $\trace\bar{D}=0$, 
\item[(b)] $\left(\n'_{\bar{X}}\bar{D}\right)\bar{Y}
-\left(\n'_{\bar{Y}}\bar{D}\right){\bar{X}}=0$
for all $\bar{X},\bar{Y}\in \mathfrak{X}(L)$ 
\end{itemize}
where
$\n'$ is the  Levi-Civita connection  of the metric induced by $g$.

Conversely, if the pair $(g,\gamma)$ is hyperbolic  (respectively, 
elliptic) with respect to a tensor $\bar{J}$ on $L^2$ satisfying 
${\bar J}^2=\bar I$ (respectively, ${\bar J}^2=-\bar I$), then the  
hypersurface  $f$ is hyperbolic (respectively, elliptic) with  respect 
to the horizontal lift $J$ of $\bar J$. In addition, the horizontal lift 
$D=\mu J$ of a tensor $\bar{D}=\bar{\mu}\bar{J}$, $\mu=\bar{\mu}\circ\pi$,  
satisfying $(a)$ and $(b)$ also fulfills the properties $(i)$--$(iv)$ in 
Lemma \ref{conditions}.
\end{lemma}

\proof We have from parts $(i)$ and $(ii)$ of Lemma \ref{conditions} 
and  Corollary $13$ in \cite{dft2} that the tensor $D$ is projectable. 
Then part $(iii)$ of Lemma \ref{conditions} gives 
$\trace\bar{D}=\trace D=0$.

From part $(iv)$ of Lemma \ref{conditions} we have that that $\det D$ 
is projectable and from Lemma~\ref{para} that also $J$ is projectable.
We have from the Gauss parametrization that
$$
f_*AX=-N_*X=-h_*\pi_*X
$$ 
where  $h=i\circ g$. Hence,
\be\label{first}
f_*ADX=h_*\pi_*DX=-h_*\bar{D}\pi_*X
\ee
for any $X\in\mathfrak{X}(M)$. In particular,
\be\label{bra}
f_*AD[X,Y]=-h_*\bar{D}\pi_*[X,Y]=-h_*\bar{D}[\pi_*X,\pi_*Y]
\ee
for any $X,Y\in\mathfrak{X}(M)$. Moreover,
\begin{eqnarray}\label{net}
f_*\n_XADY\!\!\!&=&\!\!\!\nab_Xf_*ADY-\<AX,ADY\>N\nonumber\\
\!\!\!&=&\!\!\!-\nab_{\pi_*X}h_*\bar{D}\pi_*Y
-\<h_*\pi_*X,h_*\bar{D}\pi_*Y\>h\circ\pi\nonumber\\
\!\!\!&=&\!\!\!-h_*\n'_{\pi_*X}\bar{D}\pi_*Y-\a_h(\pi_*X,\bar{D}\pi_*Y)
-\<\pi_*X,\bar{D}\pi_*Y\>h\circ\pi\nonumber\\
\!\!\!&=&\!\!\!-h_*\n'_{\pi_*X}\bar{D}\pi_*Y-\a_g(\pi_*X,\bar{D}\pi_*Y).
\end{eqnarray}
 From (\ref{6}) and the above, we have that
$$
0=f_*(\n_XB)Y-f_*(\n_YB)X=f_*\n_XADY-f_*\n_YADX-f_*AD[X,Y].
$$
We conclude that part $(b)$ holds as well as 
$$
\a_g(\pi_*X,\bar{D}\pi_*Y)=\a_g(\bar{D}\pi_*X,\pi_*Y).
$$
Since $\bar D\in\spa\{I,\bar J\}$ but $\bar D\not\in\spa\{I\}$, the preceding
equation is equivalent to
\be\label{equiv}
\a_g(\pi_*X,\bar{J}\pi_*Y)=\a_g(\bar{J}\pi_*X,\pi_*Y)
\ee
and thus $g$ is hyperbolic (respectively, elliptic) with respect to  $\bar J$.

To deal with the function $\gamma$ we first show that condition (\ref{eq:aJ})
is equivalent to
$$
(\hess h^v+h^v I)J=J^t(\hess h^v+h^v I)
$$
where $\hess h^v$ is the endomorphism  of $TL$ associated to the Hessian
and $h^v=\<h,v\>$ for any $v\in\R^{n+1}$.
We have that the Hessian of $h^v$ satisfies
$$
\hess h^v(X,Y)=\<\alpha_h(X,Y),v\>=\<i_*\alpha_g(X,Y)-\<X,Y\>h,v\>
$$
for all $X,Y\in\mathfrak{X}(L)$. Thus
$$
\<i_*\alpha_g(JX,Y)-i_*\alpha_g(X,JY),v\>
=\<((\hess {h^v}+h^v I)J-J^t(\hess {h^v}+h^v I))X,Y\>
$$
for all $X,Y\in\mathfrak{X}(L)$.

It remains to prove that
\be\label{prenda}
\left(\hess\gamma+\gamma I\right)\bar{J}
=\bar{J}^t\left(\hess\gamma+\gamma I\right).
\ee
By the Gauss parametrization, there exists a diffeomorphism 
$\theta\colon U\subset\Lambda\to M^n$ from an open neighborhood of the 
zero section of $\Lambda$ such that $\pi\circ\theta=\bar\pi$ and
$$
f\circ\theta(x,w)=\gamma(x)h(x)+h_*\n\gamma(x)+w
$$
for any $(x,w)\in\Lambda$.
Let $j\colon T_xL\to T_{(x,w)}\Lambda$ be the linear isometry in
Proposition \ref{gp}. Then,  for any $\bar{X},\bar{Y}\in T_xL$ we obtain 
using (\ref{eq}) and (\ref{first})   that
\begin{eqnarray}\label{addpw}
-\< AD\theta_* j\bar{X},\theta_*j\bar{Y}\>\!\!\!&=&\!\!\!
-\< f_ *AD\theta_* j\bar{X},f_* \theta_*j\bar{Y}\>\nonumber\\
\!\!\!&=&\!\!\!\<h_*\bar{D}\pi_* \theta_* j\bar{X},h_*\bar{Y}\>\nonumber\\
\!\!\!&=&\!\!\!\<\bar{D}\bar\pi_*j\bar{X},\bar{Y}\>'\nonumber\\
\!\!\!&=&\!\!\!\<\bar{D}P_w^{-1}\bar{X},\bar{Y}\>'.
\end{eqnarray}
It follows that $\bar{D}P_w^{-1}=P_w^{-1}\bar{D}^t$, or equivalently, that
$P_w\bar{D}=\bar{D}^tP_w$. And because $\bar D\in\spa\{I,\bar J\}$, this 
is equivalent to $P_w\bar{J}=\bar{J}^tP_w$.
Moreover, using that $A_w\bar{J}=\bar{J}^{t}A_w$ as follows from
(\ref{equiv}), we conclude that  (\ref{prenda}) is satisfied.\qed

\begin{lemma}\label{propr}\po
The following assertions on a surface  $g\colon L^2\to\Sf^n$ are equivalent:
\begin{itemize}
\item[(i)] The surface $g$ is hyperbolic (respectively, elliptic) with respect
to a tensor $\bar{J}$ on $L^2$ satisfying ${\bar J}^2=I$ (respectively,
${\bar J}^2=-I$), and there is $\bar{D}=\bar{\mu}\bar{J}$,  
$\bar\mu>0$, such that
\begin{itemize}
\item[(a)] $\trace\bar{D}=0$,
\item[(b)] $\left(\n'_{\bar{X}}\bar{D}\right)\bar{Y}
-\left(\n'_{\bar{Y}}\bar{D}\right){\bar{X}}=0$
for all $\bar{X},\bar{Y}\in \mathfrak{X}(L).$
\end{itemize}
\item[(ii)] There exist real-conjugate (respectively, complex-conjugate) 
coordinates on $L^2$ such that system  (\ref{eq:siti}) 
(respectively, (\ref{eq:comp})) has solution.
\end{itemize}
\end{lemma}

\proof  We make use of Proposition \ref{coordinates}. In the case of 
real coordinates $(u,v)$ and since we have $\bar{D}\p u=\bar\mu\p u$, 
$\bar{D}\p v=-\bar\mu\p v$, we easily see that
$$
\left(\n'_{\p u}\bar{D}\right)\p v
-\left(\n'_{\p v}\bar{D}\right)\p u=0
$$
is equivalent to the system (\ref{eq:siti}). The case of complex coordinates
is similar.\vspace{1,5ex}\qed

\noindent\emph{Proof of Theorem \ref{main}:} By Proposition \ref{hypel}, 
on each connected component of an  open and dense subset of $M^n$ the hypersurface $f$ 
is either hyperbolic or elliptic with respect to  $J\in\Gamma(\End(\Delta^\perp))$. 
It follows from Lemma \ref{conditions} that there exists $D\in\Gamma(\End(\Delta^\perp))$ 
satisfying the properties  $(i)$--$(iv)$. 

Let $f$ be parameterized by a pair $(g,\gamma)$ in terms of the Gauss 
parametrization.
By Lemma \ref{prop*} if $f$ is  hyperbolic (respectively, elliptic)
with respect to $J$, then $J$ and $D$ can be projected to  tensors $\bar{J}$ and
$\bar D\in\spa\{I,\bar{J}\}$ on $L^2$, with   $\bar{J}^2=I$ (respectively,
$\bar{J}^2=-I$).
Moreover,  the pair $(g,\gamma)$ is  hyperbolic (respectively, elliptic)  
with respect to  $\bar{J}$ and $\bar D$ satisfies $(a)$ and $(b)$ in 
Lemma~\ref{prop*}. Now the proof of the direct statement follows from 
Lemma \ref{propr}.

Conversely, let $f\colon M^n\to\R^{n+1}$ be a simply-connected 
hypersurface parameterized  in terms of the Gauss parametrization
by a special hyperbolic or special elliptic pair $(g,\gamma)$.
By Lemma~\ref{propr},  there exists $\bar{D}=\bar{\mu}\bar{J}$ satisfying 
equations $(a)$ and $(b)$. 
By Lemma~\ref{prop*} the hypersurface  $f$ is hyperbolic or elliptic 
with respect to the horizontal lift $J$ of $\bar J$, and the horizontal 
lift $D=\mu J$ of  $\bar{D}=\bar{\mu}\bar{J}$ satisfies $(i)$--$(iv)$ in  
Lemma \ref{conditions}. 

To conclude from Theorem \ref{mains} that $f$ admits a unique nontrivial 
infinitesimal bending it remains to show that $B\in\Gamma(\End(TM))$ 
defined by $B\vert_{\Delta^\perp}=A\vert_{\Delta^\perp}D$ and $\Delta\subset\ker B$
is symmetric and satisfies equations (\ref{6}) and (\ref{5}). In fact, 
that $B$ is symmetric follows easily from (\ref{addpw}). 

Part $(b)$ of Lemma \ref{prop*} together with (\ref{bra}) and (\ref{net}) 
imply that
$$
f_*((\n_XB)Y-(\n_YB)X)=f_*\n_XADY-f_*\n_YADX-f_*AD[X,Y]=0
$$
for any $X,Y\in\Gamma(\Delta^\perp)$.  
Since  $D$ is projectable, we can use Corollary $13$ in \cite{dft2} and 
deduce that
$$
\n_TD=[D,C_T]
$$
for any $T\in\Delta$. Using (\ref{codazzi}) we obtain 
$$
(\n_TB)X=(\n_TAD)X=(\n_TA)DX+A[D,C_T]X=ADC_TX
$$
for any $T\in\Delta$ and $X\in\Delta^\perp$. It follows that
$$
(\n_XB)T-(\n_TB)X=0
$$
for any $T\in\Delta$ and $X\in\Delta^\perp$. Since
$$
(\n_SB)T-(\n_TB)S=B[S,T]=0
$$
for any $S,T\in\Delta$, we have shown that (\ref{bra}) holds.
Since (\ref{net}) is equivalent to $\trace D=0$, the proof follows.\qed

\section{The ruled case}

In this section, we discuss the infinitesimal bendings of ruled 
hypersurfaces that have not been considered yet.
\vspace{1ex}

Let $f\colon M^n\to\R^{n+1}$, $n\geq 3$, be a ruled hypersurface without
flat points which is not surface-like on any open subset of $M^n$. Then $f$
has rank two and there exists locally an orthonormal frame $X,Y$ of
$\Delta^\perp$ such that the second fundamental form is the form
$$
A\vert_{\Delta^\perp}=\begin{bmatrix}
\lambda&\mu\\
\mu&0
\end{bmatrix}.
$$
Note that if $M^n$ is simply-connected then the set of all isometric 
immersions of $M^n$ into $\R^{n+1}$ consists of ruled immersions with 
the same rulings; see \cite{dft}. Moreover, this set
can be parametrized by the set of all smooth functions in an interval.
In fact, the second fundamental form of any other immersion must
be of the form
$$
A\vert_{\Delta^\perp}=\begin{bmatrix}
\lambda + \theta&\mu\\
\mu&0
\end{bmatrix}
$$
where $\theta\in C^\infty(M)$ is determined by choosing a smooth function
along an integral curve of $X$ and extending it to $M^n$ by requiring that
\be\label{codazzi2}
Y(\theta)=\<\n_XX,Y\>\theta\;\;\mbox{and}\;\;T(\theta)
=\<\n_XX,T\>\theta
\ee
for any $T\in\Delta$.

\begin{proposition}\po
 Let $f\colon M^n\to\R^{n+1}$, $n\geq 3$, be a simply-connected ruled
hypersurface of constant rank two that is not surface-like on any open 
subset of $M^n$.  Then, any infinitesimal bending is the variational 
vector field of an isometric bending.
\end{proposition}

\proof Since $M^n$ is simply-connected, there is a global orthonormal 
frame $\{X,Y\}$ of $\Delta^\perp$ as above. By Lemma \ref{conditions} and 
Proposition \ref{ruled}, the Codazzi tensor $B$ on $M^n$ is given by 
$B|_{\Delta^\perp}=A|_{\Delta^\perp}D$ and $\Delta\subset\ker B$, where 
$D=\theta J$ and  $J\in\Gamma(End(\Delta^\perp))$ is such that  $JX=Y$ and 
$JY=0$. Moreover, $\theta\in C^{\infty}(M)$ is arbitrarily prescribed along 
an integral curve of $X$ and  required to satisfy (\ref{codazzi2}). 
Therefore, the one-parameter family of Codazzi tensors $A(t)=A+tB$, $t\in\R$, 
gives rise to an isometric bending of $f$ having the infinitesimal bending 
determined by $B$ as its variational vector field.

{\renewcommand{\baselinestretch}{1}

\hspace*{-20ex}\begin{tabbing} \indent\= IMPA -- Estrada Dona Castorina, 110
\indent\indent\= Univ. of Ioannina -- Math. Dept. \\
\> 22460-320 -- Rio de Janeiro -- Brazil  \>
45110 Ioannina -- Greece \\
\> E-mail: marcos@impa.br \> E-mail: tvlachos@uoi.gr
\end{tabbing}}
\end{document}